\documentclass[10pt,a4paper]{article}
\usepackage[T1]{fontenc}
\usepackage[utf8]{inputenc}
\usepackage{amsmath,amssymb,amsfonts}
\usepackage{bm}
\usepackage{booktabs}
\usepackage{array}
\usepackage{hyperref}
\usepackage[margin=2.5cm]{geometry}

\hypersetup{
    colorlinks=true,
    linkcolor=blue,
    urlcolor=blue,
    citecolor=blue
}

\title{Structure-Driven Inversion: A New Paradigm for Solving Inverse Problems}
\author{CHEN Shengchang\\
School of Earth Sciences, Zhejiang University, Hangzhou, China\\
\texttt{chenshengc@zju.edu.cn}}
\date{}

\begin{document}

\maketitle

\begin{abstract}
Inverse problems are predominantly solved within the optimization-driven paradigm, which formulates the problem as objective-function minimization and approaches the solution by iterative search. Though universal, it suffers from limitations in efficiency, interpretability, and multi-parameter decoupling. This paper proposes a new paradigm---Structure-Driven Inversion (SDI). Instead of iterative search, SDI identifies and exploits the intrinsic structure of the problem to construct a solution method. It has two types of structure---mathematical structure and physical structure---and two corresponding driving inversion types: Mathematical Structure-Driven Inversion (MSDI) and Physical Structure-Driven Inversion (PSDI). The current representative method of MSDI is Mathematical Structure-Driven Pseudo-Inverse Inversion (MSDPII), which constructs a pseudo-inverse in the spectral domain via unitary diagonalization. The current representative methods of PSDI are Physical Structure-Driven Waveform Inversion Imaging (PSDWII), which performs inversion through virtual-source projection, and Physical Structure-Driven Back-Propagation (PSDBP), which implements three structural projections for deep neural networks via physical-system analogy. SDI is not a rejection but a complement and extension of Optimization-Driven Inversion (ODI). To the best of the author's knowledge, no existing study has systematically presented structure-driven inversion as an independent paradigm; this paper aims to fill that gap.

\end{abstract}

\section{Introduction}
\label{sec:intro}

An inverse problem is a central subject in scientific computing and engineering: the goal is to infer model parameters $m$ from observed data $d$, satisfying

\begin{equation}
F(m) = d,
\end{equation}
where $F$ is the forward operator. This is a typical ``from effects to causes'' process, and it arises in geophysical exploration, medical imaging, signal processing, system control, and many other fields. Unlike forward problems (from causes to effects), which are usually well-posed, inverse problems (from effects to causes) are intrinsically \textbf{ill-posed}: at least one of the three requirements---existence, uniqueness, and stability of the solution---fails to hold \cite{tikhonov1977}. Small perturbations in the data may be drastically amplified during inversion, rendering the solution extremely unstable.

For more than half a century, a dominant methodology has emerged for solving inverse problems---the \textbf{optimization-driven paradigm}. Its logic is to recast the inverse problem as the minimization of an objective function and to approach the optimal solution through iterative search:

\begin{equation}
\min_{m} J(m) = \|F(m) - d\|^2 + \alpha R(m),
\end{equation}
where the first term measures the least-squares data misfit and the second term $R(m)$ is a regularization term used to stabilize the solution.

The historical roots of this formulation can be traced back to Laplace and Gauss. Although they adopted different loss measures---Laplace employed the least absolute deviations criterion, while Gauss proposed the least-squares criterion---their methodological logic was identical: first construct an objective function, then minimize it to solve the inverse problem. Subsequent developments, from Tikhonov regularization to the adjoint-state method and further to SGD/Adam in deep learning, have essentially improved the objective function and the minimization algorithm within the same framework. The core logic of constructing and minimizing an objective function has remained unchanged.

The optimization-driven paradigm has dominated inverse problem solving because of its \textbf{universality}: as long as the forward operator is computable---whether differentiable or not---it can be incorporated into the framework of objective-function minimization through gradient-based methods (when differentiable) or non-gradient-based methods (such as simulated annealing and genetic algorithms, etc.). From the 1960s to the present, this paradigm has permeated every field of science and engineering. To reveal this universality, this section first summarizes the common features of optimization-driven methods across seven representative classes of inverse problems, and then analyzes their achievements and inherent limitations.

\subsection{Universality of the Optimization-Driven Paradigm: Seven Representative Classes of Inverse Problems}

The optimization-driven paradigm dominates inverse problem solving because of its cross-disciplinary universality. The following seven representative classes of inverse problems differ in their physical backgrounds and forward operators, yet they share the same methodology of ``constructing an objective function $\rightarrow$ iteratively solving''.

In \textbf{geophysical inversion}, full waveform inversion (FWI) inverts subsurface medium parameters (velocity, density, etc.) by minimizing the residual between observed and simulated waveforms, with gradients computed via the adjoint-state method \cite{tarantola1984,virieux2009}. Tarantola \cite{tarantola1984} laid the theoretical foundation of FWI and introduced the adjoint-state method to seismic inversion. Pratt et al. \cite{pratt1998} extended FWI to the frequency domain, Virieux and Operto \cite{virieux2009} provided a comprehensive review, and M\'etivier and Brossier \cite{metivier2024} gave a complete mathematical derivation of the adjoint-state method. Gravity and magnetic inversion use the solution of Poisson's equation as the forward operator, and classical solution methods include the conjugate gradient method \cite{hestenes1952} and L-BFGS \cite{liu1989}. Electromagnetic inversion recovers subsurface resistivity from surface electromagnetic data, requiring solution of Maxwell's equations \cite{haber1997}. For a systematic exposition of optimization-driven methods in geophysical inversion, the reader is referred to the classic works of Menke \cite{menke2018} and Tarantola \cite{tarantola1987,tarantola2004}.

In \textbf{image reconstruction}, computed tomography (CT) reconstruction recovers the attenuation coefficient distribution from Radon transform projection data, and total variation (TV) regularization \cite{rudin1992} effectively handles sparse-angle and low-dose problems; compressed sensing magnetic resonance imaging (MRI) achieves high-fidelity reconstruction under undersampling via $\ell_1$ regularization \cite{lustig2007}; TV-based variational methods for image deblurring and super-resolution preserve edge structures \cite{rudin1992}. For the variational methods and mathematical theory of image inverse problems, see the systematic expositions by Bertero and Boccacci \cite{bertero1998}, Scherzer et al. \cite{scherzer2009}, and Natterer \cite{natterer2001}.

In \textbf{parameter identification}, ODE/PDE parameter estimation can be formulated as a least-squares optimization problem, with core algorithms including the Levenberg-Marquardt (L-M) and Gauss-Newton methods. For the mathematical theory of inverse problems for partial differential equations, see the authoritative work of Isakov \cite{isakov2006}.

In \textbf{inverse scattering}, the Born iterative method (BIM) and the distorted Born iterative method (DBIM) \cite{habashy1993,chew1993} solve the problem through iterative linearization approximations, facing strong nonlinearity and severe ill-posedness. For a review of recent advances in inverse scattering, see the collection edited by Uhlmann \cite{uhlmann2003}.

In \textbf{optimal control}, model predictive control (MPC) solves a finite-horizon optimization problem at each time step, with gradients computed via the adjoint-state method.

In \textbf{system identification}, linear and nonlinear system identification employs subspace methods, prediction error methods (PEM), and sparse regularization \cite{ljung1999}.

In \textbf{parameter estimation in machine learning}, the training process of deep learning---constructing a loss function, computing gradients via back-propagation, and iteratively updating with stochastic gradient descent (SGD)/Adam \cite{rumelhart1986,kingma2015}---is methodologically isomorphic to classical inverse problems: forward propagation corresponds to forward modeling, and back-propagation corresponds to the adjoint-state method. For large-scale learning methods such as diffusion models and large language models, the methodological core still does not transcend the optimization-driven framework.

The common features of the above seven classes of inverse problems are summarized in Table \ref{tab:seven_problems}.

\begin{table}[htbp]
\centering
\small
\caption{Comparison of optimization-driven solutions for seven classes of inverse problems}
\label{tab:seven_problems}
\begin{tabular}{p{2.2cm}p{2cm}p{2.2cm}p{2.2cm}p{2.2cm}}
\toprule
Field & Core problem & Forward operator & Optimization-driven method & Main challenge \\
\midrule
Geophysical inversion & FWI, gravity inversion & Wave equation / Poisson equation & Adjoint-state method + gradient descent & Computational cost, multi-parameter coupling \\
Image reconstruction & CT/MRI reconstruction & Radon/Fourier sampling & TV regularization + ADMM & Large-scale, underdetermined \\
Parameter identification & ODE/PDE parameter estimation & Differential equations & L-M / Gauss-Newton & Non-identifiability, ill-posedness \\
Inverse scattering & Property reconstruction & Wave equation integral & BIM/DBIM & Strong nonlinearity, ill-posedness \\
Optimal control & MPC & State-space equations & Adjoint method + optimization solver & Real-time constraints \\
System identification & Model construction & Input-output mapping & PEM / subspace methods & High-dimensional parameters \\
Machine learning & Parameter learning & Neural network forward pass & SGD/Adam + BP & Ultra-large-scale, non-convex \\
\bottomrule
\end{tabular}
\end{table}

\subsection{Achievements and Limitations of the Optimization-Driven Paradigm}

The optimization-driven paradigm dominates all of the above fields because of its \textbf{universality}: as long as the forward operator $F$ is differentiable, one can construct an objective function, compute gradients, and solve iteratively. From Tikhonov regularization \cite{tikhonov1977} to deep learning \cite{rumelhart1986}, and from FWI to MPC, the same methodological logic runs through.

However, the limitations of the optimization-driven paradigm also stem precisely from its universality: it ignores the special structure inherent in the problem. When a translation-invariant operator is treated as a general matrix, when a physical process with a causal chain is treated as a general mapping, and when a deep network with layered organization is treated as a general composite function, efficiency, interpretability, and decoupling capability are all compromised. Moreover, in the adjoint-state method, the adjoint wavefield serves merely as a mathematical auxiliary quantity, and its physical meaning remains unclear.

From a more fundamental perspective, these limitations of the optimization-driven paradigm are not accidental. In the early twentieth century, the French mathematician Hadamard pointed out that the well-posedness of a mathematical problem depends on three criteria: existence, uniqueness, and continuous dependence of the solution on the data \cite{hadamard1923}. Forward problems---from causes to effects---usually satisfy these three criteria; inverse problems---from effects to causes---usually do not. This is the origin of the statement that ``inverse problems are ill-posed.''

Why are most inverse problems ill-posed? The root causes are as follows. First, the degrees of freedom of the observed data are far fewer than those of the model parameters, resulting in insufficient information dimensionality. Second, the forward operator is usually an integral-type or smoothing-type operator (mathematically, a compact operator), whose singular values decay rapidly, so that the division in the inversion amplifies noise without bound. Third, when the forward problem is nonlinear, the same data may correspond to multiple models. Together, these factors lead to the ill-posedness of inverse problems. In addition, inferring causes from effects lacks a corresponding mathematical and physical equation; that is, no inverse mathematical-physical equation exists that directly describes the process from effects to causes. The ``inversion equation'' for from-effects-to-causes reasoning is therefore artificially constructed.

However, inverse problems are ``mostly'' ill-posed, not ``universally'' ill-posed. When an inverse problem possesses special structure---for example, when the forward operator is a unitary transformation (information-complete and energy-conserving), or when it can be completely diagonalized by a specific transform---the inversion can be well-posed. This is precisely the starting point of the structure-driven idea in this paper: identify and exploit the intrinsic structure of the problem, so that inversion returns from ``blind search'' to ``structural construction.''

For the computational framework of optimization-driven methods, see Vogel \cite{vogel2002}; for the Bayesian framework for inverse problems, see Kaipio and Somersalo \cite{kaipio2005} (whose computational implementation also typically relies on optimization); for a comprehensive treatment of parameter estimation and inverse problems, see Aster et al. \cite{aster2018}. For the application of regularization methods in geophysics and generalized inversion, see Zhdanov \cite{zhdanov2002} and Groetsch \cite{groetsch1999}. For a unified treatment of optimization and regularization theory across various inverse problems, see Wang et al. \cite{wang2010}.

\subsection{Positioning of This Paper: A New Paradigm}

Based on the above understanding, the division of labor between the two paradigms can be articulated more precisely:

\begin{quote}
The optimization-driven paradigm, which approaches the solution through iterative search, is suitable for general ill-posed inverse problems, compensating for missing information through regularization.

The structure-driven paradigm, which constructs a solution method by exploiting the structure of the forward operator, is suitable for inverse problems with identifiable structure, using structural information to transform ill-posed inversion into well-posed construction.
\end{quote}

The relationship between the two is not substitution but complementarity: \textbf{Structure-Driven Inversion} (SDI) provides an efficient, direct path for problems with ``sufficient structure,'' while \textbf{Optimization-Driven Inversion} (ODI) maintains universal solving capability for problems with ``missing structure.''

A survey of the existing literature reveals that from Tikhonov regularization \cite{tikhonov1977} to TV regularization \cite{rudin1992}, from Bayesian inversion \cite{kaipio2005} to physics-informed neural networks (PINNs) \cite{raissi2019}, and from structurally constrained inversion \cite{haber1997,menke2018} to joint inversion---although the types of constraints employed differ ($\ell_2$ smoothness, $\ell_1$ sparsity, statistical priors, PDE residuals, geological structures, etc.)---\textbf{they share the same methodological core}:

\begin{equation}
\min_{m} \|F(m) - d\|^2 + \lambda_1 R_1(m) + \lambda_2 R_2(m) + \cdots,
\end{equation}

where the first term is the data-fidelity term and the subsequent terms are constraint terms of various types. The data-fidelity term plus the various constraint terms constitutes the objective function, which is then solved using iterative algorithms such as gradient descent, conjugate gradient, ADMM, or SGD. In a single sentence: \textbf{these methods all operate within the optimization-driven framework, augmented by mathematical and/or physical constraints}.

In contrast, SDI provides a different path. Its core is not ``define objective function $\rightarrow$ impose constraints $\rightarrow$ iterative optimization,'' but rather:

\begin{equation}
\text{identify problem structure} \rightarrow \text{exploit structure} \rightarrow \text{construct solution method} \rightarrow \text{solve}.
\end{equation}

SDI does not treat structural information as a constraint term in the objective function, but as the basis for constructing a solution method. In SDI, there are two types of structure---mathematical structure and physical structure---and two corresponding driving inversion types: \textbf{Mathematical Structure-Driven Inversion} (MSDI) and \textbf{Physical Structure-Driven Inversion} (PSDI). The current representative method of MSDI is \textbf{Mathematical Structure-Driven Pseudo-Inverse Inversion} (MSDPII; Section \ref{sec:math_structure}), which constructs a pseudo-inverse through pointwise division in the spectral domain \cite{chen2026a,chen2026b}. The current representative methods of PSDI are \textbf{Physical Structure-Driven Waveform Inversion Imaging} (PSDWII; Section \ref{sec:phys_structure}), which constructs the solution through virtual-source projection and deconvolution \cite{chen2026c}, and \textbf{Physical Structure-Driven Back-Propagation} (PSDBP; Section \ref{sec:analogy_structure}), which interprets a deep neural network as a physical system through structural analogy and constructs the solution through three structural projections \cite{chen2026d}. These three representative methods are not isolated contributions; they are the cumulative result of long-term work on structure-driven inversion, covering mathematical structure, physical structure, and analogical structure respectively, and they form the current methodological core of SDI. There is no limit to the number of specific methods under each driving inversion type; additional MSDI and PSDI methods can be developed as new structures are identified.

ODI is \textbf{universal}: it does not concern itself with the specific structure of the problem, and any differentiable forward operator can be incorporated. SDI is \textbf{universal as a methodology but customized in its methods}: the paradigm of SDI is general, but each specific method is designed for a particular structure. The relationship between the two is not substitution but complementarity: ODI is the ``master key,'' while SDI is the ``custom mold.''

To the best of the author's knowledge, no study in the literature has systematically presented ``structure-driven inversion'' as an independent methodological paradigm. This paper aims to fill that gap by providing a systematic alternative solution path for inverse problems with identifiable structure.

\section{Related Work}

Before systematically elaborating the structure-driven paradigm, this section reviews the main categories of methods related to SDI in order to clarify the boundaries of its originality.

\subsection{Regularization Methods}

Tikhonov regularization \cite{tikhonov1977} is one of the earliest methods to systematically address the ill-posedness of inverse problems; its core idea is to add a penalty term to the objective function in order to stabilize the solution. TV regularization \cite{rudin1992} further exploits the gradient sparsity of images as a prior. These methods utilize mathematical structure as a regularization term within the optimization-driven framework, but they do not change the basic logic of ``constructing an objective function $\rightarrow$ iteratively solving.''

\subsection{Statistical and Bayesian Inversion}

Bayesian inversion \cite{kaipio2005} places the inverse problem in a probabilistic framework, describing prior knowledge of the solution through a prior distribution; maximization or sampling of the posterior distribution still relies on iterative optimization or Markov chain Monte Carlo (MCMC) sampling. Its methodological foundation remains optimization- or statistical-inference-driven.

\subsection{Physics-Informed Neural Networks (PINNs)}

PINNs \cite{raissi2019} incorporate the residual of the physical equation as part of the loss function and train the network via back-propagation and gradient descent. Physical constraints act in the form of loss terms, and the methodological foundation remains optimization-driven. The ``physical information'' in PINNs functions as a regularization constraint, not as a basis for constructing a solution method.

\subsection{Structurally Constrained Inversion and Joint Inversion}

In geophysics, ``structurally constrained inversion'' usually refers to the use of geological structural information (such as horizons and faults) as constraints \cite{haber1997,menke2018}; joint inversion uses structural similarity between multiple physical fields as constraints. These methods exploit physical structural information of the solution within the optimization framework rather than stepping outside that framework.

\subsection{Model Reduction and Operator Learning}

Model reduction methods (such as proper orthogonal decomposition, POD) and operator learning methods (such as DeepONet and FNO) exploit low-dimensional structure in data or learn mappings from parameters to solutions in order to accelerate computation or replace forward modeling. They are essentially ``data-driven'' or ``model-driven'' acceleration techniques rather than changes to the inversion solution paradigm.

The common feature of all the above methods is that they exploit structural information within the optimization-driven framework, either as a regularization term, as a prior, or as an acceleration technique. Their differences lie in the type of structure exploited (mathematical, physical, statistical, or data-based), not in the fundamental methodology. \textbf{The essential difference between SDI and these methods is that the structural information concerns the mathematical or physical structure of the forward operator or system response; it is not ``encoded'' into the objective function but instead used to construct the solution method.}

\section{The Optimization-Driven Paradigm}
\label{sec:opt_paradigm}

\subsection{Unified Formulation of the Paradigm}

All methods within the optimization-driven paradigm share the same logical framework:

\begin{equation}
\text{problem} \rightarrow \text{objective function} \rightarrow \text{gradient/derivative} \rightarrow \text{iterative update} \rightarrow \text{convergence}.
\end{equation}

Although specific implementations differ, the essence is the same (Table \ref{tab:opt_methods}).

\begin{table}[htbp]
\centering
\small
\caption{Unified logic of optimization-driven methods}
\label{tab:opt_methods}
\begin{tabular}{p{3.5cm}p{4.5cm}p{4cm}}
\toprule
Method & Objective function & Update rule \\
\midrule
Tikhonov regularization \cite{tikhonov1977} & $\|Gm-d\|^2 + \alpha\|m\|^2$ & $(G^HG + \alpha I)^{-1}G^Hd$ \\
Gradient descent \cite{nocedal2006} & $J(m)$ & $m_{k+1} = m_k - \alpha \nabla J(m_k)$ \\
FWI \cite{tarantola1984,virieux2009} & $\|d_{\text{obs}} - d_{\text{syn}}(m)\|^2$ & $m_{k+1} = m_k - \alpha (\partial d/\partial m)^T \Delta d$ \\
Back-Propagation \cite{rumelhart1986} & $\text{Loss}(y_{\text{pred}}, y_{\text{true}})$ & $W_{k+1} = W_k - \eta \nabla \text{Loss}$ \\
\bottomrule
\end{tabular}
\end{table}

From Tikhonov \cite{tikhonov1977} (1960s) to Adam \cite{kingma2015} (2010s), the development represents technological progress within the same paradigm---richer algorithms, larger scales, and greater automation---but the methodological logic remains unchanged.

\subsection{Core Operation of ODI: Cross-Correlation (Multiplication)}

The core operation of ODI is \textbf{cross-correlation/inner product}:

\begin{equation}
\nabla J = \left( \frac{\partial F}{\partial m} \right)^T \cdot (F(m) - d).
\end{equation}

The essence of this operation is \textbf{multiplication}: the sensitivity is multiplied by the residual to obtain the gradient direction.

The problem is that multiplication naturally ``mixes'' the contributions of different parameters:

\begin{equation}
\label{eq:math_coupling}
\frac{\partial J}{\partial m_i} = \sum_j \left( \frac{\partial F}{\partial m_i} \right)^T \frac{\partial F}{\partial m_j} \delta m_j + \cdots.
\end{equation}

This is ``mathematical coupling''---non-physical coupling introduced by the way the gradient is computed. The remedy in optimization-driven methods is to use the inverse of the Hessian matrix to decouple the parameters, but this incurs a computational cost of $O(N^3)$.

\subsection{Three Inherent Limitations of ODI}

Based on the above discussion and analysis, the optimization-driven paradigm has three inherent limitations:

\textbf{Limitation 1: Computational efficiency bottleneck.} Each iteration requires a complete forward computation, and tens to hundreds of iterations are typically needed. For 3D FWI or high-resolution medical imaging, this computational cost is extremely high.

\textbf{Limitation 2: Ambiguous physical meaning.} The gradient $\nabla J$ is a mathematical construct; it reflects ``how the parameters should change to reduce the objective function'' rather than the direct causal relationship between the observed data and the physical process. This limits the physical interpretability of the inversion result. In particular, the physical meaning of the adjoint wavefield in the adjoint-state method used to construct the gradient remains unclear.

\textbf{Limitation 3: Multi-parameter mathematical coupling.} As shown in Eq.~\eqref{eq:math_coupling}, the cross-correlation operation in gradient computation mixes the contributions of different parameters, producing non-physical mathematical coupling. This coupling requires the inverse of the Hessian matrix for decoupling, incurring enormous computational cost.

For the mathematical foundations of ill-posedness and regularization theory for inverse problems, see \cite{engl1996,kirsch2011}.

\section{The Structure-Driven Paradigm}
\label{sec:struct_paradigm}

\subsection{Core Definition}

\textbf{Definition (Structure).} The term ``structure'' in this paper refers to the \textbf{invariance, decomposability, or hierarchy} that can be identified and exploited in the forward operator or system response of an inverse problem. It includes two types:

\begin{enumerate}
\item \textbf{Mathematical structure}: the operator can be unitarily diagonalized, as in the case of translation invariance;
\item \textbf{Physical structure}: the causal chain of the physical process can be decomposed, as in the case of wave propagation.
\end{enumerate}

A non-physical system such as a deep neural network can also be interpreted as a physical system through structural analogy, thereby yielding a physical structure in the analogical sense.

SDI is a methodology for solving inverse problems by identifying and exploiting the intrinsic structure of the problem to construct a solution method.

Its logical chain is:

\begin{equation}
\text{problem} \rightarrow \text{identify structure} \rightarrow \text{exploit structure} \rightarrow \text{construct solution method} \rightarrow \text{solution}.
\end{equation}

This logical chain characterizes the connotation of the structure-driven methodology: identify the structure of a problem, exploit it to construct a solution method, and solve the problem. SDI can be viewed as a concrete realization of this broader methodology in the context of inverse problems.

SDI is not a rejection of ODI, but a complement to and an extension of it. Both share the ultimate goal of ``fitting the data,'' but the paths to that goal differ: ODI primarily uses iterative search, whereas SDI primarily uses structural construction.

\subsection{Two Driving Inversion Types and Open Method Development}

SDI has two driving inversion types: MSDI and PSDI (Table \ref{tab:driving_forms}). Each driving inversion type can contain multiple specific methods. The methods presented in this paper---MSDPII, PSDWII, and PSDBP---are current representatives. The framework is open: additional MSDI and PSDI methods can be developed as new structures are identified.

\begin{table}[htbp]
\centering
\small
\caption{Two driving inversion types of SDI and their current representative methods}
\label{tab:driving_forms}
\begin{tabular}{p{3cm}p{3.5cm}p{4cm}p{2.2cm}}
\toprule
Driving inversion type & Structure source & Identification criterion & Current method \\
\midrule
\textbf{MSDI} & Algebraic properties of the operator & $G = T^H \Lambda T$ (unitary diagonalization or approximate unitary diagonalization) & MSDPII \\
\textbf{PSDI} & Causal chain of the physical process; or structural analogy of a non-physical system & Input $\rightarrow$ system response $\rightarrow$ output (decomposable) & PSDWII (direct); PSDBP (analogical) \\
\bottomrule
\end{tabular}
\end{table}

MSDI exploits the algebraic properties of the operator; its current representative method is MSDPII. PSDI exploits the causal chain of a physical process; its current representative methods are PSDWII, a direct application, and PSDBP, an analogical application to deep neural networks. As new structures are identified, more MSDI and PSDI methods can be developed.

\subsection{Core Operation: Division/Inversion/Deconvolution}

The core operation of SDI is \textbf{division/inversion/deconvolution} (Table \ref{tab:inverse_ops}).

\begin{table}[htbp]
\centering
\small
\caption{Inverse operations of current representative methods}
\label{tab:inverse_ops}
\begin{tabular}{p{3cm}p{2.8cm}p{2.8cm}p{4.5cm}}
\toprule
Driving inversion type & Method & Inverse operation & Mathematical form \\
\midrule
MSDI & MSDPII & Spectral-domain division & $q_\alpha(\lambda) = \dfrac{\bar{\lambda}}{|\lambda|^2 + \alpha}$ \\
PSDI (direct) & PSDWII & Deconvolution & $\delta m_i = -\tilde{V}_s \Big/ \left( \dfrac{\partial L}{\partial m_i} \cdot u_0 \right)$ \\
PSDI (analogical) & PSDBP & Function inversion & $Z_{\text{target}} = \sigma^{-1}(A_{\text{target}})$ \\
\bottomrule
\end{tabular}
\end{table}

The common effect of these operations is to separate the coupled contributions and directly extract the target parameters---achieving natural decoupling without the use of the Hessian matrix.

\subsection{Philosophical Divide Between the Two Paradigms}

Table \ref{tab:philosophy} summarizes the comparison between the two paradigms.

\begin{table}[htbp]
\centering
\small
\caption{Philosophical comparison of ODI and SDI}
\label{tab:philosophy}
\begin{tabular}{p{3.5cm}p{4.5cm}p{4.5cm}}
\toprule
Dimension & Optimization-driven & Structure-driven \\
\midrule
\textbf{Core question} & How can the output be made to approach the target? & Where does the output come from? \\
\textbf{Methodology} & Mathematical approximation in data space & Causal backtracking in structure space \\
\textbf{Essence} & Search & Construction \\
\textbf{Metaphor} & A blind man descending a mountain (repeated pathfinding) & Map navigation (direct path) \\
\bottomrule
\end{tabular}
\end{table}

\section{Mathematical Structure-Driven Pseudo-Inverse Inversion (MSDPII): Spectral-Domain Pseudo-Inverse}
\label{sec:math_structure}

The singular value decomposition (SVD) and generalized inverse theory of linear inverse problems form the foundation of MSDI methods. For a systematic treatment of SVD in inverse problems, see Menke \cite{menke2018} (Chapters 4--5) and Aster et al. \cite{aster2018} (Chapter 3); for the mathematical theory of generalized inverses, see Kirsch \cite{kirsch2011} (Chapter 2).

\subsection{Structure Identification}

Consider the linear inverse problem $Gm = d$. If $G$ can be unitarily diagonalized:

\begin{equation}
G = T^H \Lambda T, \quad \Lambda = \text{diag}(\lambda_0, \lambda_1, \ldots, \lambda_{N-1}),
\end{equation}

where $T$ is a unitary matrix ($T^H T = I$) and $\Lambda$ is a diagonal matrix. This is a broad class of matrices: all translation-invariant operators (convolution type) satisfy this property, and the discrete Fourier transform matrix is a special case \cite{chen2026a,chen2026b}.

\subsection{Structure Exploitation: Analytical SVD}

\textbf{Theorem} (Analytical SVD \cite{chen2026a}). For $G = T^H \Lambda T$, the singular value decomposition has the analytical form

\begin{equation}
U = T^H P, \quad \Sigma = |\Lambda|, \quad V^H = T,
\end{equation}

where $P = \text{diag}\{\lambda_k / |\lambda_k|\}$ is the phase diagonal matrix and $\Sigma = \text{diag}\{|\lambda_0|, \ldots, |\lambda_{N-1}|\}$.

\textbf{Proof.} Since $\Lambda = P\Sigma$, direct substitution yields $G = T^H \Lambda T = T^H P \Sigma T = U \Sigma V^H$. The unitarity of $U$ and $V$ follows from the properties of unitary matrices. $\square$

The significance of this result is that no numerical computation of the SVD is required; the complete SVD is obtained directly from the spectral decomposition.

\subsection{Structure-Driven Inverse Operator}

We construct the spectral-domain filtering factor \cite{chen2026a}:

\begin{equation}
q_\alpha(\lambda_k) = \frac{\bar{\lambda}_k}{|\lambda_k|^2 + \alpha}, \quad \alpha > 0.
\end{equation}

The spectral-domain pseudo-inverse operator is

\begin{equation}
G_\alpha^\# = T^H \text{diag}\{q_\alpha(\lambda_k)\} T.
\end{equation}

\textbf{Properties} \cite{chen2026a}:
\begin{itemize}
\item \textbf{Bounded stability}: $\|G_\alpha^\#\|_2 \leq 1/(2\sqrt{\alpha})$;
\item \textbf{Consistency}: $\lim_{\alpha \to 0^+} G_\alpha^\# = G^+$ (Moore-Penrose generalized inverse).
\end{itemize}

\textbf{Proof.} Since $T$ is unitary, $\|G_\alpha^\#\|_2 = \|Q_\alpha\|_2 = \max_k |\lambda_k|/(|\lambda_k|^2 + \alpha) \leq 1/(2\sqrt{\alpha})$. When $\lambda_k \neq 0$, $\lim_{\alpha \to 0^+} \bar{\lambda}_k/(|\lambda_k|^2 + \alpha) = 1/\lambda_k$. $\square$

\subsection{Comparison with ODI}

\begin{table}[htbp]
\centering
\small
\caption{Comparison between MSDPII and Tikhonov regularization}
\label{tab:math_vs_tikhonov}
\begin{tabular}{p{3cm}p{4.5cm}p{4.5cm}}
\toprule
 & Tikhonov regularization \cite{tikhonov1977} & MSDPII \cite{chen2026a,chen2026b} \\
\midrule
Starting point & Optimization problem & Matrix structure \\
Solution method & $(G^HG + \alpha I)^{-1}G^Hd$ & $T^H \text{diag}\{\bar{\lambda}/(|\lambda|^2+\alpha)\} T d$ \\
Complexity & $O(N^3)$ & $O(N \log N)$ (FFT implementation) \\
Iteration & Required & One-step \\
\bottomrule
\end{tabular}
\end{table}

The two approaches are numerically equivalent but methodologically different. Tikhonov regularization starts from an optimization problem, whereas MSDPII starts from the matrix structure.

\subsection{Applicability and Limitations}

MSDPII is applicable to linear inverse problems in which $G$ can be unitarily diagonalized, or can be approximately unitarily diagonalized with controllable approximation error (for example, in image-domain least-squares migration, where a locally space-invariant approximation is used to handle the spatially varying Hessian matrix \cite{chen2022}). Typical applications include:

\begin{itemize}
\item Translation-invariant systems (convolution-type inverse problems);
\item Partial differential equations with periodic boundary conditions;
\item Certain integral equations with symmetry;
\item Potential-field inverse problems such as gravity field inversion \cite{chen2026b}.
\end{itemize}

The limitation is that $G$ is required to be \textbf{square} and \textbf{unitarily diagonalizable}. For general linear systems that are non-square or non-diagonalizable, MSDPII is not applicable.

\section{Physical Structure-Driven Waveform Inversion Imaging (PSDWII): Virtual-Source Inversion}
\label{sec:phys_structure}

Gradient computation based on the adjoint-state method is the core tool of optimization-driven FWI; see the foundational work of Tarantola \cite{tarantola1987,tarantola2004}. PSDWII---virtual-source inversion---shares with the adjoint-state method the idea of ``reverse tracking of the physical process,'' but differs fundamentally in its operational path: the former approximates the inverse through structural projection, whereas the latter approaches it iteratively through gradient descent.

\subsection{Structure Identification: Causal Chain of Wave Propagation}

Wave propagation possesses a clear \textbf{physical causal chain} \cite{chen2026c}:

\begin{equation}
\begin{aligned}
& \text{source excitation} \rightarrow \text{incident wave propagation} \\
& \rightarrow \text{interaction with heterogeneity} \rightarrow \text{virtual source (with radiation pattern)} \\
& \rightarrow \text{secondary wave propagation} \rightarrow \text{reception}.
\end{aligned}
\end{equation}

This structure can be summarized as \textbf{two propagations and one virtual source}.

The key insight is that the relationship between observed data and model parameters is not a direct mapping; it is connected through the intermediate physical quantity of the ``virtual source.'' The virtual source is the product of the interaction between the incident wave and the heterogeneity---it is a \textbf{physical entity}, not a mathematical construct.

\subsection{Mathematical Formulation of the Virtual Source}

The background wavefield $u_0$ under the background model $m_0$ satisfies

\begin{equation}
L(m_0)u_0 = s.
\end{equation}

When a model perturbation $\delta m$ is present, the perturbed wavefield $\delta u$ satisfies

\begin{equation}
L(m_0)\delta u = V_s,
\end{equation}

where the virtual source is \cite{chen2026c}

\begin{equation}
V_s = - \left. \frac{\partial L}{\partial m} \right|_{m=m_0} \delta m \cdot u.
\end{equation}

Under the weak scattering approximation (Born approximation),

\begin{equation}
V_s^p = - \left. \frac{\partial L}{\partial m} \right|_{m=m_0} \delta m \cdot u_0.
\end{equation}

\subsection{Physical Classification of Virtual Sources}

According to the relationship between the heterogeneity scale and the wavelength, virtual sources can be classified into three types \cite{chen2026c}:

\begin{table}[htbp]
\centering
\small
\caption{Physical classification of virtual sources}
\label{tab:virtual_source}
\begin{tabular}{p{3.5cm}p{3.5cm}p{4.5cm}}
\toprule
Type & Scale condition & Corresponding inversion task \\
\midrule
Scattering virtual source $V_{ss}$ & Heterogeneity $\leq$ wavelength & FWI (model parameter inversion) \\
Body reflection virtual source $V_{sr}(\alpha,\theta)$ & Heterogeneity $>$ wavelength & Property imaging (impedance inversion) \\
Surface reflection virtual source $V_{sr}(r)$ & Reflection interface & Structural imaging (migration) \\
\bottomrule
\end{tabular}
\end{table}

The mathematical expressions of the three virtual sources differ, but they share the same physical structure: the incident wave interacts with the heterogeneity to excite the virtual source, and the virtual source generates the secondary wave.

\subsection{Virtual-Source Inversion}

The inversion logic of PSDWII is \textbf{reverse tracking of the physical causal chain} \cite{chen2026c}:

\begin{equation}
\text{observed data} \rightarrow \text{adjoint projection inversion of the virtual source} \rightarrow \text{model update}.
\end{equation}

The specific steps are as follows.

\textbf{Step 1}: Approximate inversion of the virtual source from the observed data:

\begin{equation}
\tilde{V}_s = G_s^* \cdot (d_{\text{obs}} - d_{\text{syn}}),
\end{equation}

where $G_s^*$ is the adjoint operator of the secondary wave propagation operator (reverse-time extrapolation), and the symbol ``$\cdot$'' denotes the action of the adjoint operator.

\textbf{Step 2}: Remove the incident wave-field and radiation patterns from the virtual source:

\begin{equation}
\delta m_i = -\tilde{V}_s \Big/ \left( \frac{\partial L}{\partial m_i} \cdot u_0 \right).
\end{equation}

In practical computation, to avoid numerical instability when the denominator approaches zero, a small parameter $\epsilon$ can be added to the denominator.

The physical meaning is that each step is a reverse tracking of the physical process: adjoint projection extracts the virtual source, and division (deconvolution) removes the radiation patterns and incident wavefield of the virtual source.

\subsection{Handling of Nonlinear Problems: Iteration Under Structural Constraints}

When the nonlinearity is strong (multiple scattering cannot be neglected), the $u$ in the virtual source contains the perturbed field $\delta u$, and inversion cannot be completed in a single step. In this case, \textbf{stepwise inversion under structural constraints} is adopted \cite{chen2026c}:

\begin{enumerate}
\item \textbf{Iteration 1}: Born approximation (ignoring $\delta u$) $\rightarrow$ $\delta m_1$ (rough estimate);
\item \textbf{Iteration 2}: Update the background model $m_1 = m_0 + \delta m_1$ $\rightarrow$ recompute the forward modeling $\rightarrow$ extract the new virtual source $\rightarrow$ $\delta m_2$ (fine correction);
\item \textbf{Iteration 3}: Repeat until the residual meets the requirement.
\end{enumerate}

The essential difference from optimization-driven iteration is shown in Table \ref{tab:iter_compare}.

\begin{table}[htbp]
\centering
\small
\caption{Comparison between optimization-driven iteration and structure-driven iteration}
\label{tab:iter_compare}
\begin{tabular}{p{3.5cm}p{4.5cm}p{4.5cm}}
\toprule
Dimension & Optimization-driven iteration & Structure-driven iteration \\
\midrule
Direction per step & Gradient descent direction & Structural projection direction \\
Physical meaning & Ambiguous & Clear (virtual-source correction) \\
Typical number of iterations & 50--200 & Fewer \\
Essence & Blind search & Structural correction \\
\bottomrule
\end{tabular}
\end{table}

\subsection{Applicability and Limitations}

PSDWII is applicable to wave-type inverse problems with a clear causal chain, including:

\begin{itemize}
\item Seismic full waveform inversion (FWI);
\item Acoustic and elastic inverse scattering problems;
\item Least-squares reverse time migration;
\item Reverse time migration (RTM).
\end{itemize}

Its limitations are as follows:

\begin{itemize}
\item It requires a relatively accurate background model (meeting kinematic accuracy; otherwise, cycle skipping occurs);
\item The validity of the Born approximation limits the deviation range between the initial model and the true model;
\item For strongly nonlinear problems, multiple iterations are still required.
\end{itemize}

\section{Physical Structure-Driven Back-Propagation (PSDBP): Three Structural Projections}
\label{sec:analogy_structure}

The second example of PSDI is its analogical application to deep neural networks. A deep neural network, although not an actual physical system, possesses a layered structure that can be understood by analogy as a physical system: inter-layer connections act as stiffness matrices, activation functions act as constitutive relations, and forward propagation is an excitation-response process. This physical-system view forms the basis of PSDBP \cite{chen2026d}, which can be regarded as a transfer of the physical-structure-driven idea of PSDWII \cite{chen2026c} to deep learning.

\subsection{Structure Identification in Deep Neural Networks}

The structure-driven idea can be applied to deep neural networks because they possess an \textbf{identifiable layered structure} \cite{chen2026d}.

\begin{table}[htbp]
\centering
\small
\caption{Structural correspondence between deep networks and physical systems \cite{chen2026d}}
\label{tab:analogy}
\begin{tabular}{p{3.5cm}p{3.5cm}p{4.5cm}}
\toprule
Physical system element & Neural network counterpart & Mathematical representation \\
\midrule
Generalized force & Inter-layer ``demand'' signal & $\Delta^{(l)}$ \\
Generalized displacement & Layer activation output & $A^{(l)}$ \\
Stiffness matrix & Inter-layer connection weights & $W^{(l)}$ \\
Constitutive relation & Activation function & $\sigma_l$ \\
External excitation & Input data & $X^{(0)}$ \\
System response & Network output & $A^{(L)}$ \\
\bottomrule
\end{tabular}
\end{table}

\subsection{Three Structural Projections}

Based on the physical structure analogy, three structural projections are constructed \cite{chen2026d}.

\textbf{Projection 1: Adjoint demand transport (geometric structure)}

\begin{equation}
\Delta^{(l-1)} = (W^{(l)})^T \Delta^{(l)}.
\end{equation}

Force transport occurs through connection geometry (weight transpose) and is independent of constitutive properties, \textbf{excluding activation derivative factors from back-propagation}.

\textbf{Projection 2: Constitutive inversion (constitutive structure)}

\begin{equation}
Z_{\text{target}}^{(l)} = \sigma_l^{-1}(A^{(l)} + \gamma_l \Delta^{(l)}).
\end{equation}

Given the target displacement, the generalized force is recovered through the analytical inverse function, \textbf{using $\sigma^{-1}$ instead of $\sigma'$}.

\textbf{Projection 3: Least-squares approximation projection (stiffness update)}

\begin{equation}
\Delta W^{(l)} = \eta \cdot E^{(l)} \cdot (X^{(l-1)})^T.
\end{equation}

The stiffness correction is determined by the input-residual correlation.

\subsection{Essential Difference from Back-Propagation}

\begin{table}[htbp]
\centering
\small
\caption{Comparison between BP and PSDBP \cite{chen2026d}}
\label{tab:bp_vs_psdbp}
\begin{tabular}{p{3cm}p{4.5cm}p{4.5cm}}
\toprule
Dimension & BP (optimization-driven) \cite{rumelhart1986} & PSDBP (structure-driven) \cite{chen2026d} \\
\midrule
Core operation & Chain-rule differentiation & Structural projection \\
Backward signal & Error derivative & Target displacement correction \\
Activation function handling & $\sigma'$ (local derivative) & $\sigma^{-1}$ (analytical inverse) \\
Backward path & Contains $\sigma'$ factor & No $\sigma'$ factor \\
Physical meaning & Ambiguous & Clear (force--displacement--constitutive) \\
\bottomrule
\end{tabular}
\end{table}

\subsection{Applicability and Limitations}

PSDBP is applicable to deep networks with layered structure, including:

\begin{itemize}
\item Fully connected feedforward networks;
\item Convolutional neural networks (CNNs);
\item Residual networks (ResNets);
\item Recurrent neural networks (RNNs).
\end{itemize}

Its limitations are as follows:

\begin{itemize}
\item The activation function is required to be strictly monotonic (invertible);
\item Stability handling of the analytical inverse of the activation function is needed during depth propagation;
\item The theoretical analysis requires further empirical verification.
\end{itemize}

\section{Unified Framework of SDI}
\label{sec:unified}

\subsection{Unified Mathematical Formulation of the Two Driving Inversion Types}

The two driving inversion types (Sections \ref{sec:math_structure}, \ref{sec:phys_structure} and \ref{sec:analogy_structure}) share the same methodological core:

\begin{equation}
\text{identify structure} \rightarrow \text{exploit structure} \rightarrow \text{construct solution method}.
\end{equation}

More specifically, the forward and inverse processes of the current representative methods can be unified as follows:

\begin{equation}
\text{Forward:} \quad m \rightarrow \text{substructure 1} \rightarrow \text{local action} \rightarrow \text{substructure 2} \rightarrow d,
\end{equation}

\begin{equation}
\text{Inverse:} \quad d \rightarrow \text{undo substructure 2} \rightarrow \text{local inversion} \rightarrow \text{undo substructure 1} \rightarrow m.
\end{equation}

The forward process proceeds along ``model $\rightarrow$ substructure 1 $\rightarrow$ local action $\rightarrow$ substructure 2 $\rightarrow$ data,'' whereas the inverse process proceeds in the opposite direction, first undoing the substructure 2 effect, then performing local inversion, and finally undoing the substructure 1 effect. The main difference between the two lies in the middle step: the forward process involves ``local action'' (multiplication), and the inverse process involves ``local inversion'' (division/deconvolution).

\begin{table}[htbp]
\centering
\small
\caption{Unified framework of the two driving inversion types and their current representative methods}
\label{tab:unified}
\begin{tabular}{p{2.8cm}p{3.5cm}p{3.5cm}p{3.5cm}}
\toprule
 & MSDI (MSDPII) & PSDI (PSDWII) & PSDI (PSDBP) \\
\midrule
Structure source & Operator algebraic properties & Physical causal chain & Structural analogy of a non-physical system \\
Substructure 1 & Forward transform (e.g., FFT) & Incident wave propagation & Forward inter-layer connection \\
Local forward operator & Multiply by $\lambda_k$ & Radiation pattern $\cdot \delta m \cdot u_0$ & Activation function $\sigma_l$ \\
Substructure 2 & Inverse transform (e.g., IFFT) & Secondary wave propagation & Activation function output \\
Local inverse operator & Divide by $\lambda_k$ (filtering) & Deconvolution by (radiation pattern $\cdot u_0$) & Inverse function $\sigma_l^{-1}$ \\
Typical application & Convolution-type inversion \cite{chen2026a,chen2026b} & Wave inversion \cite{chen2026c} & Deep learning \cite{chen2026d} \\
\bottomrule
\end{tabular}
\end{table}

\subsection{Iteration View of SDI}

Whether SDI requires iteration, and in what form, depends on the nature and completeness of the available structure. The essential difference between SDI and ODI is not ``whether to iterate'' but ``the nature of iteration'': optimization-driven iteration is a blind search, whereas structure-driven iteration is a structural projection.

\subsection{Applicability Conditions and Boundaries}

The conditions for structure-driven applicability are as follows:

\begin{enumerate}
\item \textbf{Mathematical structure}: the operator can be unitarily diagonalized (e.g., translation invariance) or approximately unitarily diagonalized;
\item \textbf{Physical structure}: the causal chain can be decomposed (e.g., wave propagation), or a non-physical system can be interpreted as a physical system through structural analogy (e.g., deep networks).
\end{enumerate}

When none of the above conditions is satisfied, SDI degenerates to ODI. The two paradigms are complementary rather than substitutive: SDI can provide efficient and transparent solutions for ``structured problems,'' while ODI maintains universal applicability for ``unstructured problems.''

\section{Comparison of SDI and ODI}

\subsection{Systematic Comparison at the Paradigm Level}

\begin{table}[htbp]
\centering
\small
\caption{Paradigm comparison between ODI and SDI}
\label{tab:paradigm_compare}
\begin{tabular}{p{3cm}p{4.5cm}p{4.5cm}}
\toprule
Dimension & Optimization-driven paradigm & Structure-driven paradigm \\
\midrule
\textbf{Core logic} & Define objective function $\rightarrow$ iterative minimization & Identify structure $\rightarrow$ construct solution method \\
\textbf{Core operation} & Cross-correlation / multiplication & Division / inversion / deconvolution \\
\textbf{Coupling handling} & Generates mathematical coupling $\rightarrow$ Hessian remedy & Only physical coupling $\rightarrow$ natural decoupling \\
\textbf{Solution method} & Iterative search (tens to hundreds) & Structural construction \\
\textbf{Physical meaning} & Ambiguous (mathematical construct) & Clear (structural process) \\
\textbf{Applicability} & Universal & Customized \\
\textbf{Theoretical basis} & Optimization theory \cite{nocedal2006} & Structural analysis + inversion theory \\
\textbf{Metaphor} & A blind man descending a mountain & Map navigation \\
\bottomrule
\end{tabular}
\end{table}

\subsection{Historical Positioning}

\begin{table}[htbp]
\centering
\footnotesize
\caption{Evolution of methodologies for inverse problem solving}
\label{tab:history}
\begin{tabular}{p{2.5cm}p{3.2cm}p{4.5cm}p{4cm}}
\toprule
Period & Paradigm & Representative method & Essence \\
\midrule
1800s & Optimization-driven (embryonic) & Least absolute deviations (Laplace); least squares (Gauss) & Construct an objective function and minimize it \\
1960s onward & Optimization-driven & Tikhonov regularization \cite{tikhonov1977} & Construct an objective function to solve \\
1980s onward & Optimization-driven & Conjugate gradient \cite{hestenes1952}, adjoint-state method \cite{tarantola1984} & Same paradigm, larger scale \\
2010s onward & Optimization-driven & SGD \cite{rumelhart1986}, Adam \cite{kingma2015}, BP & Same paradigm, more automated \\
2020s onward & \textbf{Structure-driven} & MSDPII \cite{chen2026a,chen2026b}, PSDWII \cite{chen2026c}, PSDBP \cite{chen2026d} & New paradigm \\
\bottomrule
\end{tabular}
\end{table}

From Tikhonov to Adam is technological progress within the same paradigm; from optimization-driven to structure-driven is a fundamental shift across paradigms.

\subsection{As a Paradigm Extension}

The positioning of SDI can be stated as follows:

\textbf{SDI is a complement to current ODI problem solving, and also a paradigm extension.}

\begin{itemize}
\item \textbf{As a complement}: it compensates for the shortcomings of ODI in efficiency, decoupling, and interpretability.
\item \textbf{As a paradigm extension}: it elevates the focus of inverse problem solving from ``data fitting algorithms'' to ``co-design of structure and algorithm.''
\end{itemize}

\section{Discussion: Structure-Informed Methods and the Boundary of SDI}
\label{sec:discussion}

The boundary between SDI and ODI can be further illustrated by recent structure-informed data-driven methods. A representative example is the signal-fitting diffusion denoising probabilistic model with total variation constraint (SFDDPM-TV) \cite{zhang2026sfddpm}. Starting from the wave equation and Green's function formulation, Zhang et al. \cite{zhang2026sfddpm} identify two structural properties of seismic data---structural simplicity and spatial regularity---and replace the conventional noise-fitting objective of denoising diffusion probabilistic models (DDPMs) with a signal-fitting objective, further incorporating a weighted total variation (TV) regularizer.

A closely related idea appears in image generation. Li and He \cite{li2025jit} propose Just image Transformers (JiT), which directly predicts clean images rather than noise or noised quantities, grounded in the manifold assumption: natural data lie on a low-dimensional manifold, whereas noised quantities are distributed across the full high-dimensional space \cite{li2025jit}. SFDDPM-TV independently adopts signal fitting for prestack seismic denoising, motivated by the structural simplicity and spatial regularity of seismic data derived from wave propagation, and further incorporates TV regularization to enhance cross-area generalization \cite{zhang2026sfddpm}. Both methods share the insight that fitting clean signals is fundamentally different from fitting noise, and both exploit this structural fact to redesign the learning objective.

However, these methods remain optimization-driven: their training and sampling still rely on iterative optimization and a diffusion Markov chain, and the structural prior acts as a loss-level or objective-level design rather than constructing a direct inverse operator. Their core operation is therefore not division, inversion, or deconvolution, but optimization under structural constraints. Accordingly, both SFDDPM-TV and JiT can be regarded as structure-informed optimization-driven methods rather than strict SDI methods.

This example helps clarify the complementary relation between ODI and SDI. ODI can incorporate structural priors to improve generalization and stability, whereas SDI goes further by using structure to construct the solution method itself. The two paradigms are therefore not mutually exclusive: SDI provides a direct, structure-customized path for problems with identifiable structure, while ODI remains the universal framework for problems with missing or incomplete structure. As new structure-informed methods continue to emerge, the boundary between ODI and SDI will be further tested and refined, and SDI is expected to provide a unifying perspective for distinguishing structure-constrained optimization from structure-driven construction.

\section{Conclusion}

This paper has systematically elaborated SDI, a new paradigm for inverse problem solving. The core understandings are as follows.

\begin{enumerate}
\item \textbf{SDI is a new paradigm.} It is not a replacement for ODI, but a complement to and an extension of it. Both share the ultimate goal of ``data fitting,'' but the paths to that goal differ: ODI primarily uses iterative search, whereas SDI primarily uses structural construction.

\item \textbf{SDI has two driving inversion types}: MSDI and PSDI. Each driving inversion type can contain multiple specific methods. The current representative methods are MSDPII, PSDWII, and PSDBP. MSDPII is the mathematical structure-driven pseudo-inverse inversion; PSDWII is the physical structure-driven waveform inversion imaging (direct application); and PSDBP is the physical structure-driven back-propagation (analogical application). The framework is open, and additional MSDI and PSDI methods can be developed as new structures are identified. All methods share the same methodological core: identify structure $\rightarrow$ exploit structure $\rightarrow$ construct solution method (Section \ref{sec:unified}).

\item \textbf{The core operation of SDI is division/inversion/deconvolution.} This contrasts with the core operation of ODI (cross-correlation/multiplication). Division (deconvolution) naturally decouples the contributions, avoiding the ``mathematical coupling'' problem inherent in ODI and eliminating the need for the Hessian matrix (Sections \ref{sec:opt_paradigm} and \ref{sec:struct_paradigm}).

\item \textbf{SDI may also require iteration} (Sections \ref{sec:phys_structure}, \ref{sec:analogy_structure} and \ref{sec:unified}). Whether iteration is needed, and in what form, depends on the nature and completeness of the available structure---each step being a structural projection rather than a blind search.

\item \textbf{The applicability condition of SDI is that the problem must possess an identifiable mathematical or physical structure.} When the structure is missing, the method degenerates to ODI. The two paradigms are complementary, forming a complete methodological system for inverse problem solving (Section \ref{sec:unified}).

\item \textbf{The paradigm significance of SDI}: it marks the methodological evolution of inverse problem solving---from the general framework of ``problem $\rightarrow$ objective function $\rightarrow$ optimization'' to the structured framework of ``problem $\rightarrow$ structure identification $\rightarrow$ structure-customized solution method.'' This is an important complement to and extension of the optimization-driven paradigm.
\end{enumerate}

\textbf{Final remarks.}

ODI is the ``universal search strategy'' for inverse problem solving---repeatedly finding a path in unknown terrain.

SDI is the ``customized direct strategy''---constructing a path using the problem structure.

The ultimate pursuit of SDI is to make inversion, in all possible cases, return to exploiting the intrinsic laws of the system (mathematical or physical) for direct solution, rather than approaching it through blind iteration. This is not only an improvement in efficiency but also a deepening of cognition.

Each step of the iteration has a structural basis---this is precisely the essence of SDI as a ``new paradigm.''

\end{document}